# A symplectic slice theorem

Juan-Pablo Ortega
Institut Nonlinéaire de Nice
UMR 129, CNRS-UNSA
1361, route des Lucioles
06560 Valbonne, France
Juan-Pablo.Ortega@inln.cnrs.fr

Tudor S. Ratiu
Département de Mathématiques
École Polytechnique Fédérale
CH 1015 Lausanne. Switzerland
Tudor.Ratiu@epfl.ch

October 3, 2001

**Abstract**

We provide a model for an open invariant neighborhood of any orbit in a symplectic manifold endowed with a canonical proper symmetry. Our results generalize the constructions of Marle [Mar84, Mar85] and Guillemin and Sternberg [GS84] for canonical symmetries that have an associated *momentum map*. In these papers the momentum map played a crucial role in the construction of the tubular model. The present work shows that in the construction of the tubular model it can be used the so called *Chu map* [Chu75] instead, which exists for any canonical action, unlike the momentum map. Hamilton's equations for any invariant Hamiltonian function take on a particularly simple form in these tubular variables. As an application we will find situations, that we will call *tubewise Hamiltonian*, in which the existence of a standard momentum map in invariant neighborhoods is guaranteed.

## Contents

# 1 Introduction

Let $(M, \omega)$ be a symplectic finite dimensional manifold and $G$ be a finite dimensional Lie group that acts properly and canonically on it via the left action $\Phi : G \times M \to M$. Denote by $\mathfrak{g}$ the Lie algebra of $G$ and by $\mathfrak{g}^*$ its dual. Marle [Mar84, Mar85] and Guillemin and Sternberg [GS84] have independently shown that if this action is Hamiltonian, that is, it has an associated globally defined momentum map $\mathbf{J} : M \to \mathfrak{g}^*$ then, for any point $m \in M$ there is an open $G$-invariant neighborhood of the corresponding orbit $G \cdot m$ that admits a representation in terms of a symplectic twist product in which the group $G$ appears as a factor. This representation generalizes the Slice Theorem introduced, in the category of proper $G$-spaces, by Palais [Pal61].

This letter shows that rather than the momentum map, it is the so called **Chu map** [Chu75], always available in the presence of canonical actions, which allows the statement and proof of the Symplectic Slice Theorem. These results will be presented in Section 2. Since the proof of the main result follows the classical one, once the appropriate changes have been made, our presentation will be somewhat sketchy, leaving well known facts as exercises for the reader but always providing the necessary steps to make it self contained.

In Section 3 we will write down the differential equations that describe the Hamiltonian vector field associated to any $G$-invariant Hamiltonian function, in the variables provided by the Symplectic Slice Theorem. These equations are known in the Hamiltonian context under the name of **Reconstruction equations** [O98, OR02] or **Bundle equations** [RWL99].

In the presence of a canonical action, the vector fields associated to the corresponding infinitesimal generators are always locally Hamiltonian and hence, for any point $m \in M$, there is always an open neighborhood $U$ around it where we can define a momentum map $\mathbf{J}_U : U \to \mathfrak{g}^*$ whose components are the Hamiltonian functions of the infinitesimal generator vector fields. Notice that this argument does **not** imply that every symplectic $G$–space is locally a Hamiltonian $G$–space since the neighborhood $U$ cannot always be chosen $G$–invariant (as an example take the natural action of the two-torus $\mathbb{T}^2$ on itself). In Section 4 we will use the Symplectic Slice Theorem to provide sufficient conditions under which the neighborhoods $U$ can be chosen $G$–invariant making our $G$–space *tubewise Hamiltonian*.

# 2 A normal form for canonical proper actions

Let $(M, \omega)$ be a symplectic manifold, $G$ be a Lie group acting canonically and properly on it, and $m \in M$ be an arbitrary point in $M$ around which we will construct the slice coordinates. We start the construction by stating the following readily verifiable facts:

**(i)** The vector space $V_m := T_m(G \cdot m)^\omega / (T_m(G \cdot m)^\omega \cap T_m(G \cdot m))$ is symplectic with the symplectic form $\omega_{V_m}$ defined by

$$\omega_{V_m}([v], [w]) := \omega(m)(v, w),$$

for any $[v] = \pi(v)$ and $[w] = \pi(w) \in V_m$, and where $\pi : T_m(G \cdot m)^\omega \to T_m(G \cdot m)^\omega / (T_m(G \cdot m)^\omega \cap T_m(G \cdot m))$ is the canonical projection. The vector space $V_m$ is called the ***symplectic*** ***normal space*** at $m$. The superscript $\omega$ on a vector subspace of a symplectic vector space, denotes the $\omega$-orthogonal complement.

(ii) Let $H := G_m$ be the isotropy subgroup of $m$ with Lie algebra $\mathfrak{h}$. The properness of the $G$-action guarantees that $H$ is compact. The mapping $(h, [v]) \longmapsto [h \cdot v]$, with $h \in H$ and $[v] \in V_m$, defines a linear symplectic action of the Lie group $H$ on $(V_m, \omega_{V_m})$, where $g \cdot u$ denotes the tangent lift of the $G$–action on $TM$, for $g \in G$ and $u \in TM$. We will denote by $\mathbf{J}_{V_m} : V_m \longrightarrow \mathfrak{h}^*$ the associated $H$-equivariant momentum map given by:

$$\langle \mathbf{J}_{V_m}(v), \zeta \rangle = \frac{1}{2}\omega_{V_m}(\zeta \cdot v, v) \quad \text{for all} \quad \zeta \in \mathfrak{h}.$$

(iii) The vector subspace $\mathfrak{k} \subset \mathfrak{g}$ given by

$$\mathfrak{k} = \{\eta \in \mathfrak{g} \mid \eta_M(m) \in (\mathfrak{g} \cdot m)^\omega\} \tag{2.1}$$

is a subalgebra of $\mathfrak{g}$ and $\mathfrak{h} \subset \mathfrak{k}$.

The symbol $\eta_M(m)$ denotes the evaluation at $m \in M$ of the infinitesimal generator $\eta_M$ of the action, associated to the element $\eta \in \mathfrak{g}$, that is, $\eta_M(m) = \frac{d}{dt}\big|_{t=0} \exp t\eta \cdot m$. Also, $\mathfrak{g} \cdot m := \{\eta_M(m) \mid \eta \in \mathfrak{g}\}$ is the tangent space $T_m(G \cdot m)$ at $m$ to the $G$-orbit $G \cdot m$.

The fact that $\mathfrak{k}$ is a subalgebra of $\mathfrak{g}$ is a consequence of the closedness of the symplectic form $\omega$ and the canonical character of the $G$-action. We provide below the proof.

We have to show that for any $\eta, \zeta \in \mathfrak{k}$, the bracket $[\eta, \zeta]$ also belongs to $\mathfrak{k}$, that is,

$$\omega(m)([\eta, \zeta]_M(m), \xi_M(m)) = 0 \quad \text{for all} \quad \xi \in \mathfrak{g}. \tag{2.2}$$

Recall that, in general,

$$\begin{aligned} 0 &= \mathbf{d}\omega(\eta_M, \zeta_M, \xi_M) \\ &= \eta_M[\omega(\zeta_M, \xi_M)] + \zeta_M[\omega(\xi_M, \eta_M)] + \xi_M[\omega(\eta_M, \zeta_M)] \\ &\quad -\omega([\eta_M, \zeta_M], \xi_M) + \omega([\eta_M, \xi_M], \zeta_M) - \omega([\zeta_M, \xi_M], \eta_M). \end{aligned} \tag{2.3}$$

If we evaluate the previous expression at the point $m$, the last two summands vanish since $\eta_M(m), \zeta_M(m) \in (\mathfrak{g} \cdot m)^\omega$ and $[\eta_M, \xi_M](m) = -[\eta, \xi]_M(m)$, $[\zeta_M, \xi_M](m) = -[\zeta, \xi]_M(m) \in \mathfrak{g} \cdot m$. In order to study the other terms we will need the relation that states that for any $\xi \in \mathfrak{g}$ and any $g \in G$:

$$\xi_M(g \cdot m) = \frac{d}{dt}\bigg|_{t=0} \exp t\xi \cdot (g \cdot m) = \frac{d}{dt}\bigg|_{t=0} gg^{-1} \exp t\xi \cdot (g \cdot m) = T_m\Phi_g(\mathrm{Ad}_{g^{-1}}\xi)_M(m). \tag{2.4}$$

Consequently,

$$
\begin{aligned}
(\eta_M[\omega(\zeta_M, \xi_M)])(m) &= \frac{d}{dt}\bigg|_{t=0} \omega(\exp t\eta \cdot m)\left(\zeta_M(\exp t\eta \cdot m), \xi_M(\exp t\eta \cdot m)\right) \\
&= \frac{d}{dt}\bigg|_{t=0} \omega(\exp t\eta \cdot m)\left(T_m \Phi_{\exp t\eta}\left(\mathrm{Ad}_{\exp(-t\eta)}\zeta\right)_M(m), T_m \Phi_{\exp t\eta}\left(\mathrm{Ad}_{\exp(-t\eta)}\xi\right)_M(m)\right) \\
&= \frac{d}{dt}\bigg|_{t=0} \omega(m)\left(\left(\mathrm{Ad}_{\exp(-t\eta)}\zeta\right)_M(m), \left(\mathrm{Ad}_{\exp(-t\eta)}\xi\right)_M(m)\right) \\
&= \omega(m)([\zeta,\eta]_M(m), \xi_M(m)) + \omega(m)(\zeta_M(m), [\xi,\eta]_M(m)) = \omega(m)([\zeta,\eta]_M(m), \xi_M(m)).
\end{aligned}
$$

Analogously, it is easy to check that

$$
(\zeta_M[\omega(\xi_M, \eta_M)])(m) = \omega(m)([\zeta,\eta]_M(m), \xi_M(m)) \quad \text{and} \quad (\xi_M[\omega(\eta_M, \zeta_M)])(m) = 0.
$$

All of these equalities substituted in (2.3) establish (2.2), as required. Notice that if $\eta \in \mathfrak{h}$, then $\eta_M(m) = 0$ and hence $\eta \in \mathfrak{k}$, which shows that $\mathfrak{h} \subset \mathfrak{k}$.

(iv) **The Chu map**: The Chu map is defined as the $G$-equivariant map $\Psi : M \longrightarrow Z^2(\mathfrak{g})$ (Lie algebra two cocycles) given by

$$
\Psi(m)(\xi, \eta) = \omega(m)(\xi_M(m), \eta_M(m)) \quad \text{for all} \quad \xi, \eta \in \mathfrak{g}.
$$

It is easy to verify that Chu's map is a momentum map in the sense that its level sets are preserved by the flows associated to Hamiltonian vector fields of $G$-invariant Hamiltonians. Also, it can be checked [GS84b] that

$$
\ker T_m \Psi = ([\mathfrak{g}, \mathfrak{g}] \cdot m)^\omega.
$$

(v) **Local symplectic structure for $G \times \mathfrak{k}^*$**: The compactness of the isotropy subgroup $H$ allows us to choose two $\mathrm{Ad}_H$-invariant complements: $\mathfrak{m}$ to $\mathfrak{h}$ in $\mathfrak{k}$ and $\mathfrak{q}$ to $\mathfrak{k}$ in $\mathfrak{g}$. Therefore, we have the orthogonal decompositions

$$
\mathfrak{g} = \mathfrak{h} \oplus \mathfrak{m} \oplus \mathfrak{q}, \quad \text{where} \quad \mathfrak{k} = \mathfrak{h} \oplus \mathfrak{m}, \tag{2.5}
$$

as well as their duals

$$
\mathfrak{g}^* = \mathfrak{h}^* \oplus \mathfrak{m}^* \oplus \mathfrak{q}^*, \quad \text{where} \quad \mathfrak{k}^* = \mathfrak{h}^* \oplus \mathfrak{m}^*.
$$

The previous decomposition provides an injection $i : G \times \mathfrak{k}^* \hookrightarrow G \times \mathfrak{g}^*$ that we can use to pull back the canonical symplectic form of $G \times \mathfrak{g}^* \cong T^*G$ to $G \times \mathfrak{k}^*$ in order to obtain a closed two-form $\omega_1$ on $G \times \mathfrak{k}^*$. Also, we define on $G \times \mathfrak{k}^*$ the skew-symmetric two form $\omega_2$ by

$$
\omega_2(g, \nu)\left((T_e L_g \cdot \xi, \rho), (T_e L_g \cdot \eta, \sigma)\right) = \omega(m)\left(\xi_M(m), \eta_M(m)\right) = \Psi(m)(\xi, \eta),
$$

for any $(g, \nu) \in G \times \mathfrak{k}^*$, $\xi, \eta \in \mathfrak{g}$, and $\rho, \sigma \in \mathfrak{k}^*$. It can be proved that there is a $H$-invariant neighborhood $\mathfrak{k}_r^*$ of the origin in $\mathfrak{k}^*$ such that the restriction of the form

$$\Omega := \omega_1 + \omega_2 \tag{2.6}$$

to $T_r := G \times \mathfrak{k}_r^*$ is a symplectic form for that space. Indeed, we first show that:

**$\Omega$ is closed**: the two form $\omega_1$ is clearly closed. In order to show that $\omega_2$ is closed we define for any $\xi \in \mathfrak{g}$ and $\rho \in \mathfrak{k}^*$ the vector field $(\xi, \rho) \in \mathfrak{X}(G \times \mathfrak{k}^*)$ given by $(\xi, \rho)(g, \nu) := (T_e L_g \cdot \xi, \nu)$, whose flow is $F_t(g, \nu) = (g \exp t\xi, \rho + t\nu)$. It is easy to see that the Lie bracket of two such vector fields is given by $[(\xi, \rho), (\eta, \sigma)] = ([\xi, \eta], 0)$. Using these vector fields it is easy to see that

$$\begin{aligned}
\mathbf{d}\omega_2\left((\xi, \rho), (\eta, \sigma), (\lambda, \tau)\right) &= (\xi, \rho)\left[\omega_2\left((\eta, \sigma), (\lambda, \tau)\right)\right] - (\eta, \sigma)\left[\omega_2\left((\xi, \rho), (\lambda, \tau)\right)\right] \\
&\quad + (\lambda, \tau)\left[\omega_2\left((\xi, \rho), (\eta, \sigma)\right)\right] - \omega_2([(\xi, \rho), (\eta, \sigma)], (\lambda, \tau)) \\
&\quad + \omega_2([(\xi, \rho), (\lambda, \tau)], (\eta, \sigma)) - \omega_2([(\eta, \sigma), (\lambda, \tau)], (\xi, \rho)).
\end{aligned}$$

Now notice that for any $(g, \nu) \in G \times \mathfrak{k}^*$ we have, for instance, that

$$\begin{aligned}
((\xi, \rho)\left[\omega_2\left((\eta, \sigma), (\lambda, \tau)\right)\right])(g, \nu) &= \left.\frac{d}{dt}\right|_{t=0} \omega_2(g \exp t\xi, \nu + t\rho)\left((\eta, \sigma)(g \exp t\xi, \nu + t\rho), (\lambda, \tau)(g \exp t\xi, \nu + t\rho)\right) \\
&= \left.\frac{d}{dt}\right|_{t=0} \omega(m)(\eta_M(m), \lambda_M(m)) = 0,
\end{aligned}$$

and also,

$$\omega_2([(\xi, \rho), (\eta, \sigma)], (\lambda, \tau))(g, \nu) = \omega(m)\left([\xi, \eta]_M(m), \lambda_M(m)\right).$$

The last three equalities imply that

$$\begin{aligned}
\mathbf{d}\omega_2(g, \nu)\left((T_e L_g \cdot \xi, \rho), (T_e L_g \cdot \eta, \sigma), (T_e L_g \cdot \lambda, \tau)\right) &= \\
\omega(m)\left([\xi, \eta]_M(m), \lambda_M(m)\right) + \omega(m)\left([\xi, \lambda]_M(m), \eta_M(m)\right) &- \omega(m)\left([\eta, \lambda]_M(m), \xi_M(m)\right) \\
&= \mathbf{d}\omega(m)(\xi_M(m), \eta_M(m), \lambda_M(m)) = 0,
\end{aligned}$$

which guarantees the closedness of $\omega_2$.

**$\Omega$ is locally non degenerate:** Let $\xi = \xi_1 + \xi_2$ and $\eta = \eta_1 + \eta_2$ be arbitrary elements in $\mathfrak{g}$, with $\xi_1, \eta_1 \in \mathfrak{k}$ and $\xi_2, \eta_2 \in \mathfrak{q}$. Let also $(g, \nu) \in G \times \mathfrak{k}^*$, $\rho, \sigma \in \mathfrak{k}^*$, and suppose that for all $\eta \in \mathfrak{g}$ and $\sigma \in \mathfrak{k}^*$ we have that

$$\Omega(g, \nu)\left((T_e L_g \cdot \xi, \rho), (T_e L_g \cdot \eta, \sigma)\right) = 0,$$

or equivalently:

$$\langle \sigma, \xi_1 \rangle - \langle \rho, \eta_1 \rangle + \langle \nu, [\xi, \eta] \rangle + \omega(m)\left((\xi_2)_M(m), (\eta_2)_M(m)\right) = 0. \tag{2.7}$$

We first show that when $\nu = 0$ the previous expression is non degenerate. Indeed, suppose that $\nu = 0$. Setting $\eta = 0$ in (2.7) at letting $\sigma$ vary we obtain $\xi_1 = 0$. Also, setting $\eta_2 = 0$ and letting $\eta_1$ vary we have $\rho = 0$. Finally, since $\omega(m)\left((\xi_2)_M(m), \eta_M(m)\right) = 0$ for all $\eta \in \mathfrak{g}$, $(\xi_2)_M(m) \in \mathfrak{g} \cdot m \cap (\mathfrak{g} \cdot m)^\omega$ and hence $\xi_2 \in \mathfrak{k} \cap \mathfrak{q} = \{0\}$ and, consequently, $\xi = 0$ and $\rho = 0$. Since non degeneracy is an open condition, we can choose an $\mathrm{Ad}(H)$-invariant neighborhood $\mathfrak{k}^*_r$ of zero in $\mathfrak{k}^*$ where the expression (2.7) is non degenerate. Also, as (2.7) does not depend on $G$ (the form $\Omega$ is $G$-invariant) we can conclude that expression, and consequently $\Omega$, is non degenerate on $G \times \mathfrak{k}^*_r$.

We are now in position to introduce the symplectic twist product that will constitute the slice coordinates we are after.

**Proposition 2.1 (The slice coordinates)** *Let $(M, \omega, G)$ be a symplectic $G$-manifold, $m \in M$, $H := G_m$, $V_m$ be the symplectic normal space at $m$, and $\mathfrak{m}$ be the space introduced in the splitting (2.5). Then there are $H$-invariant neighborhoods $\mathfrak{m}^*_r$ and $(V_m)_r$ of the origin in $\mathfrak{m}^*$ and $V_m$, respectively, such that the twisted product*

$$Y_r := G \times_H (\mathfrak{m}^*_r \times (V_m)_r) \tag{2.8}$$

*is a symplectic manifold acted on by the Lie group $G$ according to the expression $g \cdot [h, \eta, v] := [gh, \eta, v]$, for any $g \in G$ and any $[h, \eta, v] \in Y_r$. This action is symplectic.*

**Proof.** The symplectic form for $Y_r$ will be obtained via a standard symplectic reduction out of the symplectic forms of $T_r$ and $V_m$. First of all, consider the left action $\mathcal{R}$ of $H$ on $T_r$ given by

$$\mathcal{R}_h(g, \nu) = (gh^{-1}, \mathrm{Ad}^*_{h^{-1}}\nu), \quad h \in H, (g, \nu) \in T_r.$$

Using the definition of $\Omega$ it is straightforward to verify that this action is globally Hamiltonian on $T_r$ with equivariant momentum map $\mathbf{J}_\mathcal{R} : T_r \to \mathfrak{h}^*$, given by

$$\mathbf{J}_\mathcal{R}((g, (\eta, \rho))) = -\eta, \qquad \text{for any} \qquad (\eta, \rho) \in \mathfrak{h}^*_r \oplus \mathfrak{m}^*_r = \mathfrak{k}^*_r.$$

As we already pointed out, the $H$–action on $V_m$ is globally Hamiltonian with momentum map $\mathbf{J}_{V_m} : V_m \to \mathfrak{h}^*$. Putting together these two actions, we construct a product action of $H$ on the symplectic manifold $T_r \times V_m$, which is Hamiltonian, with $H$–equivariant momentum map $\boldsymbol{K} : T_r \times V_m \cong G \times \mathfrak{m}^*_r \times \mathfrak{h}^*_r \times V_m \to \mathfrak{h}^*$, given by the sum $\mathbf{J}_\mathcal{R} + \mathbf{J}_{V_m}$, that is,

$$\begin{array}{rcl} \boldsymbol{K}: & G \times \mathfrak{m}^*_r \times \mathfrak{h}^*_r \times V_m & \longrightarrow \quad \mathfrak{h}^* \\ & (g, \rho, \eta, v) & \longmapsto \quad \mathbf{J}_{V_m}(v) - \eta. \end{array}$$

The $H$–action on $T_r \times V_m$ is free and proper and $0 \in \mathfrak{h}^*$ is clearly a regular value of $\boldsymbol{K}$. Therefore $\boldsymbol{K}^{-1}(0)/H$ is a well–defined Marsden–Weinstein [MW74] reduced symplectic space that can be identified with $Y_r = G \times_H (\mathfrak{m}_r^* \times (V_m)_r)$ by means of the quotient diffeomorphism $L$, induced by the $H$–equivariant diffeomorphism $l$:

$$l: \begin{array}{ccc} G \times \mathfrak{m}_r^* \times (V_m)_r & \longrightarrow & \boldsymbol{K}^{-1}(0) \subset G \times \mathfrak{m}_r^* \times \mathfrak{h}_r^* \times (V_m)_r \\ (g, \rho, v) & \longmapsto & (g, \rho, \mathbf{J}_{V_m}(v), v), \end{array}$$

where the $H$-invariant neighborhood of the origin $(V_m)_r$ has been chosen so that $\mathbf{J}_{V_m}((V_m)_r) \subset \mathfrak{h}_r^*$. We define the symplectic form $\omega_{Y_r}$ on $Y_r$ as the pull back by $L$ of the reduced symplectic form $\Omega_0$ on $\boldsymbol{K}^{-1}(0)/H$. Thus, we have the following commutative diagram with the lower arrow a symplectic diffeomorphism:

$$\begin{array}{ccc} G \times \mathfrak{m}_r^* \times (V_m)_r & \xrightarrow{l} & \boldsymbol{K}^{-1}(0) \subset G \times \mathfrak{m}_r^* \times \mathfrak{h}_r^* \times (V_m)_r \\ \pi \downarrow & & \downarrow \pi_0 \\ (G \times_H (\mathfrak{m}_r^* \times (V_m)_r), \omega_{Y_r}) & \xrightarrow{L} & (\boldsymbol{K}^{-1}(0)/H, \Omega_0). \end{array} \qquad (2.9)$$

Finally, the fact that the $G$-action on $Y_r$ in the statement is symplectic is a straightforward verification. ∎

The previous proposition, together with the Constant Rank Embedding Theorem (see [Mar84, Mar85, GS84], or [SL91, O98, OR02] for a treatment similar to ours) imply in a straightforward manner that the symplectic twisted product that we just introduced can be used to model an open invariant neighborhood of any orbit in a symplectic symmetric manifold, that is, we have a Symplectic Slice Theorem.

**Theorem 2.2 (Slice Theorem for canonical proper Lie group actions)** *Let $(M, \omega)$ be a symplectic manifold and let $G$ be a Lie group acting properly and canonically on $M$. Let $m \in M$ and denote $H := G_m$. Then the manifold*

$$Y_r := G \times_H (\mathfrak{m}_r^* \times (V_m)_r) \qquad (2.10)$$

*introduced in (2.8) is a symplectic $G$–space and can be chosen such that there is a $G$–invariant neighborhood $U$ of $m$ in $M$ and an equivariant symplectomorphism $\phi : U \to Y_r$ satisfying $\phi(m) = [e, 0, 0]$.*

**Remark 2.3** The construction presented in Theorem 2.2 generalizes that of Marle [Mar84, Mar85], and Guillemin and Sternberg [GS84], in the sense that if the $G$-action on $(M, \omega)$ has an associated momentum map $\mathbf{J} : M \to \mathfrak{g}^*$, the symplectic form for $Y_r$ constructed in Proposition 2.1 using the Chu map, is identical to the one introduced in the above mentioned papers. Moreover, it is easy to check that if the momentum map $\mathbf{J}$ has $\theta : G \to \mathfrak{g}^*$ as non equivariance cocycle (recall that $\theta$ is defined by $\theta(g) := \mathbf{J}(\Phi_g(m)) - \operatorname{Ad}^*_{g^{-1}}(\mathbf{J}(m))$, for any $g \in G$, and that this definition does not depend on $m$) and $\mathbf{J}(m) = \mu$, $G \cdot m$ being the orbit around which we construct the tubular coordinates, then the canonical $G$-action on $Y_r$ is actually Hamiltonian with momentum map $\mathbf{J}_{Y_r} : Y_r \to \mathfrak{g}^*$ given by

$$\mathbf{J}_{Y_r}([g,\eta,v]) = \operatorname{Ad}^*_{g^{-1}}(\mu + \eta + \mathbf{J}_{V_m}(v)) + \theta(g), \quad \text{for all} \quad [g,\eta,v] \in Y_r;$$

this coincides with the expression provided by [Mar84, Mar85], as expected. ♦

**Remark 2.4 (Normal form for Abelian quasi-Hamiltonian $G$-spaces)** The slice theorem presented in Theorem 2.2 can be used to provide a normal form for the Lie group valued momentum maps associated to Abelian symplectic actions. The Lie group valued momentum maps have been introduced by McDuff [McD88], in the case of circle actions, and by Alekseev, *et al* [AMM97] in the general case, in order to deal with some canonical symmetries for which it is impossible to define a standard momentum map with values in the dual of a Lie algebra. In the non Abelian case their definition is not compatible with the symplecticity of $(M, \omega)$ so we will leave aside this situation for future work. In the Abelian case, the $G$-valued momentum map $\mathbf{K} : M \to G$ associated to a symplectic $G$-action on $(M, \omega)$ is defined as the map such that for any $\xi \in \mathfrak{g}$ it satisfies

$$\mathbf{i}_{\xi_M}\omega = \mathbf{K}^*(\lambda, \xi),$$

where $(\cdot, \cdot)$ is an invariant positive definite inner product on $\mathfrak{g}$ and $\lambda \in \Omega^1(G, \mathfrak{g})$ is the invariant Maurer-Cartan form.

Assume that this map exists for the $G$-action on $(M, \omega)$ and that for a fixed element $m \in M$ we have that $\mathbf{K}(m) = g_0 \in G$. An easy exercise proves that in this case the $G$-action on the tubular model $Y_r$ constructed around the orbit $G \cdot m$ has also a $G$-valued momentum map whose expression is

$$\mathbf{K}([g,\eta,v]) = g_0 \cdot g \cdot \exp \eta \cdot \exp \mathbf{J}_{V_m}(v), \quad \text{for all} \quad [g,\eta,v] \in Y_r. \; ♦$$

## 3 The reconstruction equations for canonical proper Lie group actions

The reconstruction equations are the differential equations that determine the Hamiltonian vector field associated to a $G$-invariant Hamiltonian in the coordinates provided by Theorem 2.2. In the globally Hamiltonian context these equations can be found in [O98, RWL99, OR02]. As we will see, it is remarkable that in the absence of a momentum map, the reconstruction equations written using the symplectic form of (2.8) are *formally* identical to the ones obtained for the globally Hamiltonian case.

In order to present the reconstruction equations let $h \in C^\infty(Y_r)^G$ be the Hamiltonian function whose associated vector field $X_h$ we want to write down. Let $\pi : G \times \mathfrak{m}_r^* \times (V_m)_r \to G \times_H (\mathfrak{m}_r^* \times (V_m)_r)$ be the canonical projection. The $G$-invariance of $h$ implies that $h \circ \pi \in$ $C^\infty(G \times \mathfrak{m}_r^* \times (V_m)_r)^H$ can be understood as a $H$–invariant function that depends only on the $\mathfrak{m}_r^*$ and $(V_m)_r$ variables, that is,

$$h \circ \pi \in C^\infty(\mathfrak{m}_r^* \times (V_m)_r)^H.$$

Now, since the projection $\pi$ is a surjective submersion, the Hamiltonian vector field $X_h$ can be locally expressed as

$$X_h = T\pi(X_G, X_{\mathfrak{m}^*}, X_{V_m}),$$

with $X_G$, $X_{\mathfrak{m}^*}$, and $X_{V_m}$ locally defined smooth maps on $Y_r$ and having values in $TG$, $T\mathfrak{m}_r^*$ and $T(V_m)_r$ respectively. Moreover, using the $\mathrm{Ad}_H$–invariant decomposition of the Lie algebra $\mathfrak{g}$ introduced in the previous section, the map $X_G$ can be written, for any $[g, \rho, v] \in Y_r$, as

$$X_G([g, \rho, v]) = T_e L_g\big(X_\mathfrak{h}([g, \rho, v]) + X_\mathfrak{m}([g, \rho, v]) + X_\mathfrak{q}([g, \rho, v])\big),$$

with $X_\mathfrak{h}$, $X_\mathfrak{m}$, and $X_\mathfrak{q}$, locally defined smooth maps on $Y_r$ with values in $\mathfrak{h}$, $\mathfrak{m}$, and $\mathfrak{q}$, respectively. In what follows we give the expressions that determine $X_G$, $X_{\mathfrak{m}^*}$, and $X_{V_m}$ as a function of the differential of the Hamiltonian $h$.

First, the construction of $\mathfrak{q}$ as the orthogonal complement to $\mathfrak{k}$ guarantees that the bilinear pairing $\langle \cdot, \cdot \rangle_\mathfrak{q}$ in $\mathfrak{q}$ defined using the Chu map by

$$\langle \xi, \eta \rangle_\mathfrak{q} := \Psi(m)(\xi, \eta) = \omega(m)\left(\xi_M(m), \eta_M(m)\right)$$

is non degenerate. Let $\mathbb{P}_{\mathfrak{h}^*}$, $\mathbb{P}_{\mathfrak{m}^*}$, and $\mathbb{P}_{\mathfrak{q}^*}$ denote the projections from $\mathfrak{g}^*$ onto $\mathfrak{h}^*$, $\mathfrak{m}^*$, and $\mathfrak{q}^*$, respectively, according to the $\mathrm{Ad}_H^*$–invariant splitting $\mathfrak{g}^* = \mathfrak{h}^* \oplus \mathfrak{m}^* \oplus \mathfrak{q}^*$. The non degeneracy of $\langle \cdot, \cdot \rangle_\mathfrak{q}$ implies that the mapping

$$\begin{aligned} F : \mathfrak{k} \times \mathfrak{k}^* \times \mathfrak{q} &\longrightarrow \mathfrak{q}^* \\ (\xi, \lambda, \tau) &\longmapsto \mathbb{P}_{\mathfrak{q}^*}\left(\mathrm{ad}^*_{(\xi+\tau)}\lambda\right) + \langle \tau, \cdot \rangle_\mathfrak{q}, \end{aligned} \quad (3.1)$$

is such that $F(0, 0, 0) = 0$ and that its derivative $D_\mathfrak{q} F(0, 0, 0) : \mathfrak{q}^* \to \mathfrak{q}^*$ is a linear isomorphism. The Implicit Function Theorem implies the existence of a locally defined function $\tau : \mathfrak{k} \times \mathfrak{k}^* \to \mathfrak{q}$ around the origin such that $\tau(0, 0) = 0$, and $\mathbb{P}_{\mathfrak{q}^*}\left(\mathrm{ad}^*_{(\xi+\tau(\xi,\lambda))}\lambda\right) + \langle \tau(\xi, \lambda), \cdot \rangle_\mathfrak{q} = 0$. Let now $\psi : \mathfrak{m}^* \times V_m \to \mathfrak{q}$ be the locally defined function given by:

$$\psi(\rho, v) := \tau\left(D_{\mathfrak{m}^*}(h \circ \pi)(\rho, v), \rho + \mathbf{J}_{V_m}(v)\right).$$

Using these expressions it can be easily verified that $X_h$ is given by

$$X_G([g, \rho, v]) = T_e L_g \left(\psi(\rho, v) + D_{\mathfrak{m}^*}(h \circ \pi)(\rho, v)\right) \quad (3.2)$$

$$X_{\mathfrak{m}^*}([g, \rho, v]) = \mathbb{P}_{\mathfrak{m}^*}\left(\mathrm{ad}^*_{D_{\mathfrak{m}^*}(h \circ \pi)}\rho\right) + \mathrm{ad}^*_{D_{\mathfrak{m}^*}(h \circ \pi)}\mathbf{J}_{V_m}(v) + \mathbb{P}_{\mathfrak{m}^*}\left(\mathrm{ad}^*_{\psi(\rho,v)}(\rho + \mathbf{J}_{V_m}(v))\right) \quad (3.3)$$

$$X_{V_m}([g, \rho, v]) = B^\sharp_{V_m}(D_{V_m}(h \circ \pi)(\rho, v))) \quad (3.4)$$

The symbol $B^\sharp_{V_m} : T^* V_m \to T V_m$ denotes the vector bundle map associated to the symplectic form $\omega_{V_m}$ in $V_m$.

# 4 Tubewise Hamiltonian actions

In order to explain what we are after in this section we start with the following definition.

**Definition 4.1** *Let $(M, \omega)$ be a symplectic manifold acted canonically upon by a Lie group $G$. For any point $m \in M$, we say that the $G$–action on $M$ is **tubewise Hamiltonian** at $m$ if there exists a $G$–invariant open neighborhood of the orbit $G \cdot m$ such that the restriction of the action to the symplectic manifold $(U, \omega|_U)$ is Hamiltonian (it has an associated globally defined momentum map).*

In this section we will use the Symplectic Slice Theorem to provide sufficient conditions that ensure that a given action is tubewise Hamiltonian. This is of much use in the study of singular dual pairs [O01].

We start by recalling recall that in the proper actions framework, Theorem 2.2 guarantees that any orbit of a symplectic $G$–space, $(M, \omega)$ has an invariant neighborhood around it that can be modeled by an associated bundle like the one presented in (2.10). Consequently, we can conclude that the canonical proper $G$–action on $(M, \omega)$ is tubewise Hamiltonian if the $G$–action on each $G$–invariant model (2.10) is Hamiltonian. The following result provides a sufficient condition for that to happen.

**Proposition 4.2** *Let $(M, \omega)$ be a symplectic manifold and let $G$ be a Lie group with Lie algebra $\mathfrak{g}$, acting properly and canonically on $M$. Let $m \in M$, $H := G_m$ and $Y_r := G \times_H (\mathfrak{m}_r^* \times (V_m)_r)$ be the slice model around the orbit $G \cdot m$ introduced in (2.10). If the $G$–equivariant, $\mathfrak{g}^*$–valued one form $\gamma \in \Omega^1(G; \mathfrak{g}^*)$ defined by*

$$\langle \gamma(g) \cdot T_e L_g \cdot \eta, \xi \rangle := -\omega(m) \left( \left( \mathrm{Ad}_{g^{-1}} \xi \right)_M (m), \eta_M(m) \right) \quad \text{for any} \quad g \in G, \, \xi, \eta \in \mathfrak{g} \qquad (4.1)$$

*is exact, then the $G$–action on $Y_r$ given by $g \cdot [h, \eta, v] := [gh, \eta, v]$, for any $g \in G$ and any $[h, \eta, v] \in Y_r$, has a standard associated momentum map and thus the $G$–action on $(M, \omega)$ is tubewise Hamiltonian at $m$.*

**Proof.** By looking at the details of the construction of the symplectic form of $Y_r$ in Proposition 2.1, it is easy to see that the existence of a standard momentum map for the $G$–action on $Y_r$ is guaranteed by the existence of a momentum map for the $G$–action on the symplectic manifold $(G \times \mathfrak{k}_r^*, \Omega|_{G \times \mathfrak{k}_r^*})$ introduced in (2.6). This action is given by $g \cdot (h, \eta) := (gh, \eta)$, for any $g, h \in G$, $\eta \in \mathfrak{k}^*$. The existence of this momentum map is in turn equivalent to the vanishing of the map (see for instance page 18 of [W76])

$$[\xi] \in \mathfrak{g}/[\mathfrak{g}, \mathfrak{g}] \longmapsto \left[ \mathbf{i}_{\xi_{G \times \mathfrak{k}^*}} \Omega \right] \in H^1(G \times \mathfrak{k}^*), \quad \text{for any} \quad \xi \in \mathfrak{g}. \qquad (4.2)$$

It is easy to see that for any $\xi, \eta \in \mathfrak{g}$, $g \in G$, and $\nu, \sigma \in \mathfrak{g}^*$

$$\mathbf{i}_{\xi_{G \times \mathfrak{k}^*}} \Omega(g, \nu) \left( T_e L_g \cdot \eta, \sigma \right) = \langle \sigma, \mathrm{Ad}_{g^{-1}} \xi \rangle + \langle \nu, \left[ \mathrm{Ad}_{g^{-1}} \xi, \eta \right] \rangle + \omega(m) \left( \left( \mathrm{Ad}_{g^{-1}} \xi \right)_M (m), \eta_M(m) \right).$$

The first two terms on the right hand side of the previous expression are the differential of the real function $f \in C^\infty(G \times \mathfrak{k}^*)$ given by

$$f(g, \nu) := \langle \nu, \mathrm{Ad}_{g^{-1}} \xi \rangle,$$

hence the vanishing of (4.2) is equivalent to the exactness of the $\mathfrak{g}^*$-valued one-form $\gamma$ in the statement. ∎

The following proposition provides a characterization of the exactness of (4.1) and therefore gives another sufficient condition for the tubewise Hamiltonian character of the action.

**Proposition 4.3** *Assume the hypotheses of Proposition 4.2. Let $m \in M$, $H := G_m$ and $Y_r := G \times_H (\mathfrak{m}_r^* \times (V_m)_r)$ be the slice model around the orbit $G \cdot m$ introduced in (2.10). Let $\Sigma : \mathfrak{g} \times \mathfrak{g} \longrightarrow \mathbb{R}$ be the two cocycle induced by the Chu map, that is:*

$$\Sigma(\xi, \eta) = \omega(m) \left( \xi_M(m), \eta_M(m) \right), \quad \xi, \eta \in \mathfrak{g},$$

*and let $\Sigma^\flat : \mathfrak{g} \to \mathfrak{g}^*$ be defined by $\Sigma^\flat(\xi) = \Sigma(\xi, \cdot)$, $\xi \in \mathfrak{g}$. Then the form (4.1) is exact if and only if there exists a $\mathfrak{g}^*$-valued group one cocycle $\theta : G \to \mathfrak{g}^*$ such that*

$$T_e \theta = \Sigma^\flat.$$

*In such case the action is tubewise Hamiltonian at the point $m$. Also, in the presence of this cocycle, the map $\mathbf{J}_\theta : G \times \mathfrak{k}^* \to \mathfrak{g}^*$ given by*

$$\mathbf{J}_\theta(g, \nu) := \mathrm{Ad}^*_{g^{-1}} \nu - \theta(g) \tag{4.3}$$

*is a momentum map for the $G$–action on the presymplectic manifold $G \times \mathfrak{k}^*$ with non equivariance cocycle equal to $-\theta$.*

**Proof.** Suppose first that the form $\gamma$ in (4.1) is exact. In such case, there exists a function $\theta : G \to \mathfrak{g}^*$ such that

$$\gamma(g) = \mathbf{d}\theta(g), \tag{4.4}$$

that is, for any $\xi, \eta \in \mathfrak{g}$ and $g \in G$ we have that

$$\langle T_g \theta \cdot T_e L_g \cdot \eta, \xi \rangle = -\omega(m) \left( \left(\mathrm{Ad}_{g^{-1}} \xi\right)_M (m), \eta_M(m) \right). \tag{4.5}$$

This expression determines uniquely the derivative of $\theta$ and hence choosing $\theta(e) = 0$ fixes the map $\theta : G \to \mathfrak{g}^*$. We now show that $\theta$ is a cocycle by checking that it satisfies the cocycle identity. Indeed, for any $g, h \in G$ and any $\xi, \eta \in \mathfrak{g}$, we have

$$\begin{aligned}
\langle T_g \left(\theta \circ L_h\right) \cdot T_e L_g \cdot \eta, \xi \rangle &= \langle T_{hg} \theta \cdot T_e L_{hg} \cdot \eta, \xi \rangle = \langle \gamma(hg) \cdot T_e L_{hg} \cdot \eta, \xi \rangle \\
&= \langle \mathrm{Ad}^*_{h^{-1}} \left(\gamma(g) \cdot T_e L_g \cdot \eta\right), \xi \rangle \\
&= \langle T_g \left(\mathrm{Ad}^*_{h^{-1}} \circ \theta\right) \cdot T_e L_g \cdot \eta, \xi \rangle.
\end{aligned}$$

Therefore, for any $g, h \in G$ we have $T_g \left( \theta \circ L_h \right) = T_g \left( \operatorname{Ad}^*_{h^{-1}} \circ \theta \right)$ and consequently

$$\theta \circ L_h = \operatorname{Ad}^*_{h^{-1}} \circ \theta + c(h, n),$$

where $c(h, n) \in \mathfrak{g}^*$, for any $h \in G$, and any $n \in [1, \operatorname{Card}(G/G^\circ)]$. Or, equivalently, for any $g, h \in G$ we can write

$$\theta(hg) = \operatorname{Ad}^*_{h^{-1}} \circ \theta(g) + c(h, n). \tag{4.6}$$

If we make in this equality $g = e$ and we use that $\theta(e) = 0$ we obtain that $\theta(h) = c(h, n)$, for all $h \in G$ and $n \in [1, \operatorname{Card}(G/G^\circ)]$ and hence (4.6) becomes

$$\theta(hg) = \operatorname{Ad}^*_{h^{-1}} \circ \theta(g) + \theta(h),$$

as required. Finally, notice that from (4.5) it is easy to see that $T_e \theta = \gamma(e) = \Sigma^\flat$ and therefore $\theta$ is the one cocycle in the statement of the proposition. The converse is straightforward.

The fact that the expression (4.3) produces a momentum map for the $G$–action on $G \times \mathfrak{k}^*$ follows as a straightforward verification of the equality

$$\mathbf{i}_{\xi_{G \times \mathfrak{k}^*}} \Omega = \mathbf{d} \langle \mathbf{J}_\theta, \xi \rangle, \quad \text{for any} \quad \xi \in \mathfrak{g}. \quad \blacksquare$$

The following corollary presents two situations in which the hypotheses of Proposition 4.2 are trivially satisfied.

**Corollary 4.4** *Let $(M, \omega)$ be a symplectic manifold and let $G$ be a Lie group with Lie algebra $\mathfrak{g}$, acting properly and canonically on $M$. If either,*

**(i)** $H^1(G) = 0$, *or*

**(ii)** *the orbit $G \cdot m$ is isotropic*

*then, the $G$–action on $(M, \omega)$ is tubewise Hamiltonian at $m$.*

**Acknowledgments** We are grateful for useful and interesting conversations with James Montaldi. This research was partially supported by the European Commission and the Swiss Federal Government through funding for the Research Training Network *Mechanics and Symmetry in Europe* (MASIE). Tudor Ratiu also acknowledges the partial support of the Swiss National Science Foundation.